\documentclass{amsart}
\usepackage{amsmath}
\usepackage{txfonts}
\usepackage{amsfonts}
\usepackage{hyperref}

\newtheorem{theorem}{Theorem}[section]

\newtheorem{corollary}[theorem]{Corollary}

\theoremstyle{definition}

\newtheorem{example}[theorem]{Example}

\theoremstyle{remark}
\newtheorem{remark}[theorem]{Remark}

\numberwithin{equation}{section}

\begin{document}

\title{A result on  Ricci curvature and the second Betti number }

\author{}
\address{ Department of Mathematics, Northwest University, Xi'an 710127, China
}
\curraddr{} \email{wanj\_m@aliyun.com; wan@nwu.edu.cn}
\thanks{The research is supported by the National Natural Science Foundation of China
N0.11301416}

\author{Jianming Wan}
\address{}
\email{}
\thanks{}


\date{}

\dedicatory{}


\begin{abstract}
We prove that the second Betti number of a compact Riemannian manifold vanishes under certain Ricci curved restriction. As consequences we obtain an interesting curved restriction for compact K\"{a}hler-Einstein manifolds and a homology sphere theorem in $\dim=4, 5$.
\end{abstract}

\maketitle

\section{introduction}
The study of relation between curvature and topology is the central topic in Riemannian geometry. One of the strong tool is Bochner technique. It plays a very important role in understanding relation between curvature and Betti numbers. The first result in this field is Bochner's classical result (c.f. \cite{[W]})

\begin{theorem}
(Bochner 1946) Let $M$ be a compact Riemannian manifold with Ricci curvature $Ric_{M}>0$. Then the first Betti number $b_{1}(M)=0$.
\end{theorem}

Berger investigated that in what case the second Betti number vanishes. He proved the following (c.f. \cite{[B1]}, also see \cite{[BS]} theorem 2.8)

\begin{theorem}
(Berger) Let $M$ be a compact Riemannian manifold of
dimension $n\geq5$. Suppose that $n$ is odd and the sectional curvature satisfies that $\frac{n-3}{4n-9}\leq K_{M}<1$.
 Then the second Betti number $b_{2}(M)=0$.
\end{theorem}

Consider a different curvature condition, Micallef and Wang proved (c.f. \cite{[MW]}, also see \cite{[BS]} theorem 2.7)

\begin{theorem}
(Micallef-Wang) Let $M$ be a compact Riemannian
manifold of dimension $n\geq4$. Suppose that n is even and $M$ has positive
isotropic curvature. Then the second Betti number $b_{2}(M)=0$.
\end{theorem}

Here positive isotropic curvature means, for any four othonormal vectors $e_{1}, e_{2},e_{3},e_{4}\in T_{p}M$ , the curvature tensor satisfies
$$R_{1313}+R_{1414}+R_{2323}+R_{2424}>2|R_{1234}|.$$

Recall that the Rauch-Berger-Klingenberg's sphere theorem (c.f. \cite{[B1]}) states that a simple connected compact Riemannian manifold is homeomorphic to a sphere if the sectional curvatures lie in $(\frac{1}{4}, 1]$. A generalization of sphere theorem (dues to Micallef-Moore c.f. \cite{[MM]}) says that a compact simply
connected Riemannian manifold with positive isotropic curvature is a homotopy sphere. Hence with the help of Poincare conjecture it is homeomorphic to a sphere. From the two theorems we know that theorems 1.2 and 1.3 can not cover too many examples. In this note we shall use Ricci curvature to give a relaxedly sufficient condition for the second Betti number vanishing. Our main result is

\begin{theorem}
Let $M$ be a compact Riemannian manifold. The dimension $\dim(M)=2m$ or $2m+1$. Let $\bar{k}$ (resp. $\underline{k}$) be the maximal (resp. minimal)
sectional curvature of $M$. If the Ricci curvature of $M$ satisfies that
\begin{equation}
Ric_{M}>\bar{k}+\frac{2m-2}{3}(\bar{k}-\underline{k}),
\end{equation}
then the second Betti number $b_{2}(M)=0$.
\end{theorem}

Particularly, if $M$ is a compact Riemannian manifold with nonnegative sectional curvature, then the second Betti number vanishes provided
\begin{equation}
Ric_{M}>\frac{2m+1}{3}\bar{k}.
\end{equation}

Note that there is no dimensional restriction in theorem 1.4.

Any compact K\"{a}hler manifold does not satisfy (1.1) since it has $b_{2}\geq 1$.

The condition 1.1 is a Ricci pinching condition.  We mention that several other Ricci pinching type theorems obtained by Gu and Xu (c.f. \cite{[GX]} \cite{[XG]}, ).

As an immediate consequence, we obtain a curvature restriction for special Einstein manifolds.

\begin{corollary}
Let $M$ be a compact Einstein manifold with nonzero second Betti number. Then the Ricci curvature satisfies
\begin{equation}
Ric\leq\bar{k}+\frac{2m-2}{3}(\bar{k}-\underline{k}).
\end{equation}
In addition, if the sectional curvature is nonnegative, one must have
\begin{equation}
Ric\leq\frac{2m+1}{3}\bar{k}.
\end{equation}
Particularly (1.3) holds for any compact K\"{a}hler-Einstein manifold.
\end{corollary}

\begin{remark}
1) The condition (1.1) implies that the maximal sectional curvature $\bar{k}>0$: If $\bar{k}\leq0$, then $$\bar{k}\geq Ric_{M}>\bar{k}+\frac{2m-2}{3}(\bar{k}-\underline{k}).$$ We get $\bar{k}<\underline{k}$. This is a contradiction.

2) Since $\bar{k}>0$, of course (1.1) implies $Ric_{M}>0$.

3) If the minimal sectional curvature $\underline{k}<0$. Since $\bar{k}>0$. If $\dim(M)=2m+1$, from $$2m\bar{k}\geq Ric_{M}>\bar{k}+\frac{2m-2}{3}(\bar{k}-\underline{k}),$$ one has
$$\bar{k}>\frac{2m-2}{4m-1}|\underline{k}|.$$ Similarly $$\bar{k}>\frac{1}{2}|\underline{k}|$$ provided $\dim(M)=2m$.
\end{remark}

We use theorem 1.4 to test some simple examples.

\begin{example}
1) The space form $S^{n}$, $\bar{k}=\underline{k}=1$, $Ric=n-1=\bar{k}$ for $n=2$ and $Ric=n-1>\bar{k}$ for $n\neq2$, $b_{2}(S^{2})=1$ and $b_{2}(S^{n})=0$ for $n\neq2$.

2) $S^{2}\times S^{2}$ with product metric, $\bar{k}=1,\underline{k}=0$, $Ric=1<\bar{k}+\frac{2n-2}{3}(\bar{k}-\underline{k})$, $b_{2}(S^{2}\times S^{2})=2$.

3) $S^{m}\times S^{m}, m>4$ with product metric, $\bar{k}=1,\underline{k}=0$, $Ric=m-1>\frac{2m+1}{3}\bar{k}$, $b_{2}=0$.

4) $\mathbb{C}\mathbb{P}^{n}$ with Fubini-Study metric, $\bar{k}=4,\underline{k}=1$, $Ric=2n+2=\bar{k}+\frac{2n-2}{3}(\bar{k}-\underline{k})$, $b_{2}(\mathbb{C}\mathbb{P}^{n})=1$.

\end{example}

From the examples we know that the inequality (1.1) is sharp.

The proof of theorem 1.4 is also based on Bochner technique. But comparing with Berger and Micallef-Wang's results, we consider a different side. This allows us get a uniform result (without dimensional restriction).

\section{Proof of the theorem}
\subsection{Bochner formula} Let $M$ be a compact Riemannian manifold. Let $$\Delta=d\delta+\delta d$$ be the Hodge-Laplacian, where $d$ is the exterior differentiation and $\delta$
is the adjoint to $d$.

Let $\varphi\in\Omega^{k}(M)$ be a smooth $k$-form. Then we have the well-known Weitzenb\"{o}ck formula (c.f. \cite{[W]})
\begin{equation}
\Delta\varphi=\sum_{i}\nabla^{2}_{v_{i}v_{i}}\varphi-\sum_{i,j}\omega^{i}\wedge i(v_{j})R_{v_{i}v_{j}}\varphi,
\end{equation}
here $\nabla^{2}_{XY}=\nabla_{X}\nabla_{Y}-\nabla_{\nabla_{X}Y}$ and $R_{XY}=-\nabla_{X}\nabla_{Y}+\nabla_{Y}\nabla_{X}+\nabla_{[X,Y]}$. The $\{v_{i},1\leq i\leq n\}$ are the local orthonormal vector fields and $\{\omega_{i},1\leq i\leq n\}$ are the duality.

A $k$-form $\varphi$ is called harmonic if $\Delta \varphi=0$.

The famous \emph{Hodge theorem} states that the de Rham cohomology $H^{k}_{dR}(M)$ is isomorphic to the space spanned by $k$-harmonic forms.

Let $\varphi=\sum_{i,j}\varphi_{ij}\omega^{i}\wedge\omega^{j}$ be a harmonic 2-form. By (2.1), under the normal frame we can get (c.f. \cite{[BS]} or \cite{[B1]})
 \begin{equation}
\Delta\varphi_{ij}=\sum_{k}(Ric_{ik}\varphi_{kj}+Ric_{jk}\varphi_{ik})-2\sum_{k,l}R_{ikjl}\varphi_{kl},
\end{equation}
where $R_{ijkl}=\langle R(v_{i},v_{j})v_{k},v_{l}\rangle$ is the curvature tensor and $Ric_{ij}=\sum_{k}\langle R(v_{k},v_{i})v_{k},v_{j}\rangle$ is the Ricci tensor.

So we have
\begin{eqnarray*}
\Delta|\varphi|^{2}& = & 2\sum_{i,j}\varphi_{ij}\Delta\varphi_{ij}+2\sum_{i,j}\sum_{k}(v_{k}\varphi_{ij})^{2}\\
                   &\geq & 2\sum_{i,j}\varphi_{ij}\Delta\varphi_{ij}\\
                    &\triangleq& 2F(\varphi).\\
\end{eqnarray*}

Note that by (2.1) one has the global form of above formula $$0=-\langle\Delta\varphi,\varphi\rangle=\sum_{i}|\nabla_{v_{i}}\varphi|^{2}+\langle\sum_{i,j}\omega^{i}\wedge i(v_{j})R_{v_{i}v_{j}}\varphi,\varphi\rangle-\frac{1}{2}\Delta|\varphi|^{2}.$$  The $F(\varphi)$ is just the curvature term $\langle\sum_{i,j}\omega^{i}\wedge i(v_{j})R_{v_{i}v_{j}}\varphi,\varphi\rangle$.

\subsection{Proof of Theorem 1.4}

By Hodge theorem, we only need to show that every harmonic 2-form vanishes.

\textbf{Case 1:} Assume $\dim(M)=2m$. For any $p\in M$, we can choose an orthonormal basis $\{v_{1}, w_{1},...,v_{m},w_{m}\}$ of $T_{p}M$ such that $\varphi(p)=\sum_{\alpha}\lambda_{\alpha}v_{\alpha}^{*}\wedge w_{\alpha}^{*}$ (for instance c.f. \cite{[B1]} or \cite{[BS]}). Here $\{v_{\alpha}^{*}, w_{\alpha}^{*}\}$ is the dual basis.
Then
\begin{equation}
F(\varphi)= \sum_{\alpha=1}^{m}\lambda_{\alpha}^{2}[Ric(v_{\alpha},v_{\alpha})+Ric(w_{\alpha},w_{\alpha})]-2\sum_{\alpha,\beta=1}^{m}\lambda_{\alpha}\lambda_{\beta}
R(v_{\alpha},w_{\alpha},v_{\beta},w_{\beta})
\end{equation}

The term
\begin{eqnarray*}
-2\sum_{\alpha,\beta=1}^{m}\lambda_{\alpha}\lambda_{\beta}
R(v_{\alpha},w_{\alpha},v_{\beta},w_{\beta})
           & = &-2\sum_{\alpha\neq\beta}\lambda_{\alpha}\cdot\lambda_{\beta}\cdot
           R(v_{\alpha},w_{\alpha},v_{\beta},w_{\beta})-2\sum_{\alpha=1}^{m}\lambda_{\alpha}^{2}R(v_{\alpha},w_{\alpha},v_{\alpha},w_{\alpha})\\
           & \geq & -\frac{4}{3}(\bar{k}-\underline{k})
           \sum_{\alpha\neq\beta}|\lambda_{\alpha}|\cdot|\lambda_{\beta}|-2\bar{k}\sum_{\alpha=1}^{m}\lambda_{\alpha}^{2}\\
           & \geq &-\frac{2}{3}(\bar{k}-\underline{k})
           \sum_{\alpha\neq\beta}(\lambda_{\alpha}^{2}+\lambda_{\beta}^{2})-2\bar{k}|\varphi|^{2}\\
           & = & -\frac{2}{3}(\bar{k}-\underline{k})(2m-2)|\varphi|^{2}
           -2\bar{k}|\varphi|^{2}\\
           & = &-2[\bar{k}+\frac{2m-2}{3}(\bar{k}-\underline{k})]|\varphi|^{2}.
\end{eqnarray*}
The first $"\geq"$ follows from Berger's inequality (c.f. \cite{[B1]}): For any orthonormal 4-frames $\{e_{1},e_{2},e_{3},e_{4}\}$, one has $$|R(e_{1},e_{2},e_{3},e_{4})|\leq\frac{2}{3}(\bar{k}-\underline{k}).$$

On the other hand, by the condition (1.1) we have $$\sum_{\alpha=1}^{m}\lambda_{\alpha}^{2}[Ric(v_{\alpha},v_{\alpha})+Ric(w_{\alpha},w_{\alpha})]\geq2[\bar{k}+\frac{2m-2}{3}(\bar{k}-\underline{k})]|\varphi|^{2},$$
the equality holds if and only if $\varphi(p)=0$.

This leads to
$$F(\varphi)\geq0$$ with equality if and only if $\varphi(p)=0$. Since $$\int_{M}F(\varphi)\leq\frac{1}{4}\int_{M}\Delta|\varphi|^{2}=0,$$ we get $$F(\varphi)\equiv0.$$ Thus the harmonic 2-form $\varphi\equiv0$.

\textbf{Case 2:} If $\dim(M)=2m+1$.  For any $p\in M$, we also can choose an orthonormal basis $\{u, v_{1}, w_{1},...,v_{m},w_{m}\}$ of $T_{p}M$ such that $\varphi(p)=\sum_{\alpha}\lambda_{\alpha}v_{\alpha}^{*}\wedge w_{\alpha}^{*}$ (c.f. \cite{[B1]} or \cite{[BS]}). We also have  $$F(\varphi)= \sum_{\alpha=1}^{m}\lambda_{\alpha}^{2}[Ric(v_{\alpha},v_{\alpha})+Ric(w_{\alpha},w_{\alpha})]-2\sum_{\alpha,\beta=1}^{m}\lambda_{\alpha}\lambda_{\beta}
R(v_{\alpha},w_{\alpha},v_{\beta},w_{\beta}).$$ Thus the argument is same to the even dimensional case.

This completes the proof of the theorem.

\section{Sphere theorem in $\dim 4$ and 5}

\begin{theorem}
Let $M$ be a compact Riemannian manifold. $\dim M=4$ or 5. If $$Ric_{M}>\frac{5\bar{k}-2\underline{k}}{3},$$
then $M$ is a real homology sphere, i.e. $b_{i}(M)=0$ for $1\leq i \leq \dim M-1$.
\end{theorem}

\begin{proof}
Since $Ric_{M}>0$, from theorem 1.1 we know that $b_{1}(M)=0$. Theorem 1.4 implies that $b_{2}(M)=0$. With the help of Poincare duality, we obtain the theorem.
\end{proof}

Finally we metion a differential sphere theorem for Ricci curvature obtained by Gu and Xu ( c.f. \cite{[GX]} theorem D).

\begin{theorem}
Let $M$ be a simple conncted compact Riemannian $n$-manifold. If $$Ric_{M}>(n-\frac{11}{5})\bar{k},$$ then $M$ is diffeomorphic to $S^{n}$.
\end{theorem}

\bibliographystyle{amsplain}

\end{document}